\documentclass[12pt,letterpaper]{article}

\usepackage{amssymb}
\usepackage{amsmath}
\usepackage{amsthm}
\usepackage{amscd}
\usepackage{xypic}
\theoremstyle{plain}
\newtheorem{theorem}{Theorem}[section]
\newtheorem{proposition}{Proposition}[section]
\newtheorem{corollary}{Corollary}[section]
\newtheorem{lemma}{Lemma}[section]
\theoremstyle{remark}
\newtheorem{remark}{Remark}[section]
\theoremstyle{definition}
\newtheorem{definition}{Definition}[section]

 \newcommand{\refpaper}[3]{{\sc #1,} \emph{#2,}  #3.}

\numberwithin{equation}{section}

\renewcommand{\baselinestretch}{1.2}

\def\qed{\vrule height 5pt width 5pt depth 0pt}

\def\varinjlim_#1{\lim\limits_{\longrightarrow\atop{#1}}}

\def\st{\mathop{\rm Stab}\nolimits}
\def\rist{\mathop{\rm rist}\nolimits}
\def\mod{\mathop{\rm mod}\nolimits}
\def\Aut{\mathop{\rm Aut}\nolimits}

\def\deg{\mathop{\rm deg}\nolimits}

\def\Ker{\mathop{\rm Ker}\nolimits}
\def\Sym{\mathop{\rm Sym}\nolimits}

\begin{document}

\author{Ekaterina Pervova\footnote{This research is supported by INTAS YSF
03-55-1423. The author would like to acknowledge also support from
the Swiss National Fund for scientific research and support from
the Chair of Discrete Algebraic Mathematics (EPFL).}}
\title{Profinite completions of some groups acting on trees}
\date{}
\maketitle
{\renewcommand{\baselinestretch}{1.0}
\begin{abstract}
We investigate the profinite completions of a certain family of
groups acting on trees. It turns out that for some of the groups
considered, the completions coincide with the closures of the
groups in the full group of tree automorphisms. However, we
introduce an infinite series of groups for which that is not so,
and describe the kernels of natural homomorphisms of the profinite
completions onto the aforementioned closures of respective groups.
\end{abstract}}

\tableofcontents

\section*{Introduction}

A profinite group is {\em just infinite} if its every proper
(continuous) quotient is finite; equivalently, if every closed
normal subgroup is open. It is hereditarily just infinite if every
open subgroup is just infinite. It is known, further, that any
finitely generated profinite group which is virtually a pro-$p$
group can be mapped onto a virtually pro-$p$ just infinite group
\cite{Gr00}.

Another result from \cite{Gr00}, based partly on the results of
\cite{Wil71}, states that any profinite just infinite
group either contains an open normal subgroup which is isomorphic to the
direct product of a finite number of copies of some hereditarily just
infinite profinite group, or is a profinite branch group. The latter groups
can be defined as profinite groups with a tree-like structure lattice of
subnormal subgroups \cite{Wil00}.

It is such groups that are the main subject of the present paper.
More specifically, we start with certain just infinite {\em
self-similar groups} which possess a branching structure in the
sense mentioned above (a tree-like structure lattice of subnormal
subgroups). Such a group $G$ can, in particular, be interpreted as
a group acting on a regular rooted tree $T$.

At the start, we do not require $G$ to be profinite. Nevertheless,
as a subset of $\Aut T$, the full automorphism group of $T$, $G$
is equipped with the induced topology $\tau_1$. It is proven in
\cite{Gr00} that the closure in $\Aut T$ (i.e. the completion with
respect to $\tau_1$) of any branch group is a profinite branch
group. However, $G$ has its own profinite topology $\tau_2$. Then
a very natural question to ask, and it is the first question
considered in the paper, appears to be, {\em do the two topologies
coincide (for a given $G$)?}

This question admits an equivalent form. It is easy to see that
the level stabilizers of vertices of $T$ form a base of profinite
topology in $\Aut T$ (i.e. a system of neighborhoods of the
identity). Define {\em principal congruence subgroups} in $G$ to
be the intersections of those level stabilizers with $G$. Let us
say that $G$ has {\em congruence property} if every finite index
subgroup of $G$ is a congruence subgroup, i.e. if it contains some
principal congruence subgroup. The two topologies coincide if and
only if $G$ has congruence property \cite{BG02}, \cite{GrigProf}.
In the paper we establish that groups from a rather well-known
class, $GGS$-groups \cite{GS83, Bau93}, do possess this property.
This is done in Section~\ref{ggs_completions}.

The main part of the paper, however, is devoted to studying
profinite completions of a certain class of self-similar
$p$-groups without the congruence property. The existence of such
groups is somewhat unexpected and answers negatively Question 3 of
\cite{BGS02}. Once their existence is proven, however, the next
interesting question is to investigate the natural homomorphism of
the profinite completion of such a group onto its closure in $\Aut
T$. Namely, we describe the kernel of this homomorphism, which
turns out to be a profinite abelian group of prime exponent.
Section~\ref{egs_completions} is devoted to these questions.

The main results of the paper are stated in
Theorems~\ref{cong_nonsymggs}, \ref{no_cong_property},
\ref{EGScompletion}, and \ref{final_result}. The key steps of the
prolonged proof of Theorem~\ref{final_result} are
Lemma~\ref{first_char}, Corollaries~\ref{form_of_sequences} and
\ref{description}, and Proposition~\ref{unique_sequence}. All
necessary definitions are given in
Section~\ref{defin}.\vspace{3mm}

\noindent \textbf{Acknowledgements } Many people greatly helped me
during the work on this paper. I'd like to express my gratitude to
all of them, particularly to Rostislav Grigorchuk, who first
attracted my attention to this set of questions, and to Laurent
Bartholdi, Pierre de la Harpe, Said Sidki, and Vladimir Tarkaev
for numerous and inspiring discussions on these and related
subjects. A substantial part of this work was carried out during
my visit to University of Geneve, whose staff I would like to
thank for welcoming atmosphere.

\section{Main definitions}\label{defin}

A self-similar group is, roughly speaking, a group acting in a
self-replicating manner on the set of all words in a finite
alphabet. More precisely, choose a finite set $A$. Denote by $A^*$
the set of all words of finite length in alphabet $A$.

\begin{definition}\label{self-similar}
A group $G$ acting on $A^*$ is called {\em self-similar} if for
every $a\in A$ and $g\in G$ there exist $b\in A$ and $h\in G$ such
that
$$(aw)^g=b(w^h),$$
no matter what $w\in A^*$ we take.
\end{definition}

For historical reasons, we prefer to speak about regular rooted
trees and their vertices rather than about words. Indeed, the set
of all words in a finite alphabet can be naturally identified with
a rooted spherically homogeneous tree, where the words correspond
to vertices, the root is the empty word, and two vertices are
joined by an edge if and only if they have the form $a_1a_2\ldots
a_n$ and $a_1a_2\ldots a_na_{n+1}$ for some $n$ and some $a_i\in
A$. The number $n$ is called the {\em length} of a vertex
$u=a_1a_2\ldots a_n$ and is denoted by $|u|$. The set of all
vertices of length $n$ is called the $n$th level of $T$.

Suppose that $u=\hat{a}_1\hat{a}_2\ldots \hat{a}_n$ is a vertex.
The set of all vertices of the form
$\hat{a}_1\hat{a}_2\ldots\hat{a}_na_{n+1}a_{n+2}\ldots a_{n+m}$,
where $m\in N$ and $a_{n+i}$ range over the set $A$, forms a
subtree of $T$. We will denote that subtree by $T_u$. It is easy
to see that $T_u$ is naturally isomorphic to the same tree $T$ via
the map
$$\hat{a}_1\hat{a}_2\ldots\hat{a}_na_{n+1}a_{n+2}\ldots a_{n+m}\mapsto
a_{n+1}a_{n+2}\ldots a_{n+m}.$$ This map allows to identify
subtrees $T_v$ for all vertices $v$, with one fixed tree $T$.

Consider now an arbitrary subgroup $G$ in $\Aut T$ and a vertex
$v$ of $T$. The {\em stabilizer} of $v$ in $G$ is the subgroup
$$\st_G(v)=\{g\in G\,|\,v^g=v\}.$$
Denote also by $\st_G(n)$ the subgroup $\cap_{|v|=n}\st_G(v)$, which
keeps all vertices of level $n$ fixed.

Subgroups $\st_G(n)$ are called {\em principal congruence subgroups}
in $G$. A subgroup of $G$ which contains a principal congruence
subgroup is called a {\em congruence subgroup}.

\begin{definition}
A group $G\leq \Aut T$ is said to possess the {\em congruence
property} if any its finite index subgroup is a congruence
subgroup.
\end{definition}

It is easy to see that $\Aut T$ admits a natural map $\phi:\Aut
T\rightarrow \Aut T\wr\Sym(A)$, where $\Sym(A)$ is the group of
all permutations of elements of $A$. Thus, every element $x$ of
$\Aut T$ is given by an element $f_x\in\underbrace{\Aut T\times
\ldots \times \Aut T}_{|A|}$ and a permutation $\pi_x\in\Sym(A)$.
The latter permutation is called the {\em accompanying
permutation}, or the {\em activity}, of $x$ at the root. We write
that
$$\phi(x)=f_x\cdot\pi_x.$$

In particular, the restriction of $\phi$ onto $\st_{\Aut T}(1)$ is
an embedding (in fact, an isomorphism) of $\st_{\Aut T}(1)$ into
the direct product of $|A|$ copies of $\Aut T$. We will denote
this restriction by $\Phi_1$.

It is evident now that $\Phi_1(\st_{\Aut
T}(2))=\underbrace{\st_{\Aut T}(1)\times \ldots \times \st_{\Aut
T}(1)}_{|A|}$. Hence we can obtain an isomorphism
$$\Phi_2=(\underbrace{\Phi_1\times \ldots \times \Phi_1}_{|A|})\circ\Phi_1:
\st_{\Aut T}(2)\rightarrow\underbrace{\Aut T\times \ldots \times
\Aut T}_{|A|^2}.$$ Proceeding in this manner, we define for each
positive integer $n$ an isomorphism
$$\Phi_n=(\underbrace{\Phi_{n-1}\times \ldots \times \Phi_{n-1}}_{|A|})\circ\Phi_1:
\st_{\Aut T}(n)\rightarrow\underbrace{\Aut T\times \ldots \times
\Aut T}_{|A|^n}. $$

The above notations allow us to introduce several modifications of
the notion of self-similarity for groups acting on trees.

\begin{definition}\label{recursive}
A group $G\leq\Aut T$ is called {\em recursive} if $\phi(G)$ is
contained in $G\wr\Sym(A)$ and the map $G\rightarrow
G\wr\Sym(A)\rightarrow\Sym(A)$ is onto a transitive subgroup of
$\Sym(A)$ (the latter map is the projection).
\end{definition}

Geometrically speaking, the latter condition means that $G$ acts
transitively on each level of $T$.

\begin{definition}\label{weakly_recurrent}
A group $G\leq\Aut T$ is called {\em weakly recurrent} if it is
recursive and the set-map $G\rightarrow G\wr\Sym(A)\rightarrow G$
is onto for each coordinate.
\end{definition}

\begin{definition}\label{s_self-similar}
A group $G\leq\Aut T$ is called {\em recurrent} if it is recursive
and $\Phi_1(\st_G(1))$ is subdirect product of $|A|$ copies of
$G$, i.e. if $\st_G(1)\rightarrow \underbrace{G\times\ldots\times
G}_{|A|}\rightarrow G$ is onto for each coordinate (the latter map
is the projection).
\end{definition}

The proposition below readily follows from the above definition:

\begin{proposition}\label{recurrent}
Let $G\leq\Aut T$ be a recurrent group. Then for every $n$
$\Phi_n(\st_G(n))$ is a subdirect product of $|A|^n$ copies of
$G$.
\end{proposition}

It remains to introduce a few more notations before proceeding to
the subject matter. Notice that to every $x\in\st_{\Aut T}(v)$ we
can assign a unique automorphism $x_v\in\Aut T_v$ which is
obtained by taking the restriction of $x$ onto the subtree $T_v$
(the notation $x@v$ is also used for $x_v$). Moreover, for each
vertex $v$ $x_v$ can be considered as an element of $\Aut T$.
Thus, for each vertex $v$ there is a fixed homomorphism
$\varphi_v:\st_{\Aut T}(v)\rightarrow \Aut T$. We will denote
$\varphi_v(\st_G(v))$ by $G_v$.

In general, for different vertices $u$, $v$ of the same length
$G_u$ and $G_v$ are conjugate in $\Aut T$ but may not coincide.
For recursive groups they are conjugate in $G$, and for recurrent
groups they coincide with $G$ (this readily follows from
Proposition~\ref{recurrent}). The latter case holds for all
examples considered in this paper.

The {\em rigid stabilizer} of $v$ in $G$ is the subgroup
$$\rist_G(v)=\{\,g\in G\,|\mbox{ for any }u\in T\setminus T_v\;u^g=u\}.$$
We also denote by $\rist_G(n)$ the subgroup
$\prod_{|v|=n}\rist_G(n)$. This is, of course, a normal subgroup
in $G$ (unlike the rigid stabilizer of just one vertex). Given an
element $g\in\rist_G(v)_v$, we often denote by $g*v$ the element
of $\rist_G(v)$ such that $(g*v)_v=g$.

\begin{definition}\label{branch}
A group $G\leq\Aut T$ is called {\em branch} if for all $n$ the
index $|G:\rist_G(n)|$ is finite.
\end{definition}

Indeed, the tree-like structure lattice of subnormal subgroups
mentioned in the Introduction, is given by various $\rist_G(v)$.

Let us now describe our major examples.

\subsection{Periodic $GGS$-groups}\label{defin_ggs}

The term $GGS$-group was introduced by Gilbert Baumslag
\cite{Bau93} and refers to Rostislav {\bf\em G}rigorchuk, Narain
{\bf\em G}upta, and Said {\bf\em S}idki. Those groups act on
regular $p$-trees, where $p\geq 3$. We will always assume that $p$
is prime.

$GGS$-groups are a partial case of a wider class of groups
introduced in \cite{GS84}. Each $GGS$-group is fully defined by a
vector
\begin{equation}\label{accom_vector}
\bar{\alpha}=(\alpha_1,\alpha_2,\ldots,\alpha_{p-1})\in
\underbrace{\mathbb{GF}(p)\oplus\mathbb{GF}(p)\oplus\ldots\oplus\mathbb{GF}(p)}_{p-1}
\end{equation}
in the following way.

Construct the regular $p$-tree $T_p$ as the tree over the sequence
$(A,A,\ldots)$, $A=\{0,1,2,\ldots,p-1\}$. Each $GGS$-group
$G_{\bar{\alpha}}$ acting on $T_p$ is generated by two
automorphisms $a$ and $b$. The decomposition of $a$ is given by
$$\phi(a)=1\cdot\pi,$$
where $\pi$ is the cyclic permutation $(01\ldots p-1)$. The second
generator $b$ is in the stabilizer of level 1, and
$$\Phi_1(b)=
(a^{\alpha_1},a^{\alpha_2},\ldots,a^{\alpha_{p-1}},b).$$

The vector (\ref{accom_vector}) is called the {\em accompanying
vector} of $G_{\bar{\alpha}}$. The accompanying vector is called
symmetric if $\alpha_i=\alpha_{p-i}$ for all
$i=1,2,\ldots,\frac{p-1}{2}$. (Throughout the paper we will only
deal with nonsymmetric accompanying vectors.)

It is immediately seen that each $G_{\bar{\alpha}}$ is recurrent.
(Indeed, since it is evidently recursive, it is sufficient to
verify conditions of Definition~\ref{s_self-similar} for one
coordinate only.) In particular, it implies that for each vertex
$v$ $(G_{\bar{\alpha}})_v=G_{\bar{\alpha}}$. By \cite{Vov00}, all
groups $G_{\bar{\alpha}}$ are just infinite; as follows from
\cite{GS84} and \cite{Vov00}, they are periodic if and only if
$\sum_{i=1}^{p-1}\alpha_i=0$.

By Proposition 1 from \cite{Per02}, we have

\begin{lemma}\label{commutant_for_ggs}
Let $G_{\bar{\alpha}}$ be a $GGS$-group. Then
$$G_{\bar{\alpha}}/[G_{\bar{\alpha}},G_{\bar{\alpha}}]\cong
\langle a\rangle\times \langle b\rangle\cong \mathbb{Z}_p\times
\mathbb{Z}_p.$$
\end{lemma}

Since knowing the structure of rigid stabilizers is often crucial
for studying groups acting on trees, we cite the following
proposition proved in \cite{congruence}.

\begin{proposition}\label{cost_ggs}
If the accompanying vector of a periodic $GGS$-group
$G_{\bar{\alpha}}$ is nonsymmetric, then for any vertex $u$
$\rist_{G_{\bar{\alpha}}}(u)_u$ contains
$[G_{\bar{\alpha}},G_{\bar{\alpha}}]$.
\end{proposition}

We do not reproduce the full proof here, but it is important to
notice that the proof is based on the lemma below, which is
obtained as a partial case of Theorem~2.7.1 from \cite{RosDoc}.

\begin{lemma}\label{costt_ggs}
Let $G_{\bar{\alpha}}$ be a periodic $GGS$-group with nonsymmetric
accompanying vector $\bar{\alpha}=(\alpha_1,\ldots,\alpha_{p-1})$.
Then we have inclusion
\begin{equation}\label{rist_ggs}
\Phi_1(\gamma_3(G_{\bar{\alpha}}))\geq
\underbrace{[G_{\bar{\alpha}},G_{\bar{\alpha}}]\times\ldots\times
[G_{\bar{\alpha}},G_{\bar{\alpha}}]}_p.
\end{equation}
\end{lemma}

Proposition~\ref{cost_ggs} and the property of $G_{\bar{\alpha}}$
to be recurrent imply that $G_{\bar{\alpha}}$ is branch.

Later on we will need the following length function on $G_{\bar{\alpha}}$.
Each element $x\in G_{\bar{\alpha}}$ can be represented (not uniquely) by
a word of the form,
$$a^{\delta_1}b^{\beta_1}\ldots a^{\delta_m}b^{\beta_m}a^{\delta_{m+1}},$$
where $\delta_2,\ldots,\delta_m,\beta_1,\ldots,\beta_m\in
\mathbb{GF}(p)\setminus\{0\}$ and
$\delta_1,\delta_{m+1}\in\mathbb{GF}(p)$. The number $m$ will be
called the length of such a word. The length
$\partial_{G_{\bar{\alpha}}}(x)$ is defined as the minimal length
of words representing $x$.

Also, denote by $\partial_a(x)$ the sum of all $\delta_i$ and by
$\partial_b(x)$ the sum of all $\beta_i$ in some representation of
$x$ (both sums are considered as elements of $\mathbb{GF}(p)$). By
Lemma~\ref{commutant_for_ggs}, the sums are independent of the
choice of a particular representation, so $\partial_a(x)$ and
$\partial_b(x)$ are well-defined.

\subsection{$EGS$-groups}\label{defin_egs}

The term $EGS$-group stands for ``extended Gupta-Sidki group'', but the
word ``extended'' should not be understood in the usual algebraic sense.
While each $EGS$-group contains some specific $GGS$-group as a
subgroup, it is not an extension of it. However, in a certain sense
a periodic $GGS$-group, or, more precisely, its accompanying vector, does
define a unique $EGS$-group.

Let $G_{\bar{\alpha}}$ be a periodic $GGS$-group with generators
$a$ and $b$ and
accompanying vector $(\alpha_1,\alpha_2,\ldots,\alpha_{p-1})$. The
corresponding $EGS$-group $\Gamma_{\bar{\alpha}}$ is generated by
automorphisms $a$, $b$, and the automorphism $c$ such that
$$\Phi_1(c)=(c,a^{\alpha_1},a^{\alpha_2},\ldots,a^{\alpha_{p-1}}).$$
It is immediately obvious that each $\Gamma_{\bar{\alpha}}$ is
recurrent. In particular, it implies that for each vertex $v$
$(\Gamma_{\bar{\alpha}})_v=\Gamma_{\bar{\alpha}}$.

The group $G_{\bar{\alpha}}$ is called the {\em associated
$GGS$-group} of $\Gamma_{\bar{\alpha}}$. We also denote by
$F_{\bar{\alpha}}$ the subgroup generated by $a$ and $c$ and prove
the following lemma:

\begin{lemma}\label{conjugate}
Subgroups $G_{\bar{\alpha}}$ and $F_{\bar{\alpha}}$ are conjugate
in $\Aut T$.
\end{lemma}

\begin{proof}
It is easy to see that they are conjugated by the automorphism
$\mathcal{C}=af$, where $f$ is in $\st_{\Aut T}(1)$ and is defined
by the equality,
$$\Phi_1(f)=(\underbrace{\mathcal{C},\mathcal{C},\ldots,\mathcal{C}}_p).$$
Indeed, we have that
$$a^{\mathcal{C}}=a,$$
$$b^{\mathcal{C}}=(b^{\mathcal{C}},(a^{\alpha_1})^{\mathcal{C}},\ldots,
(a^{\alpha_{p-1}})^{\mathcal{C}})=(b^{\mathcal{C}},a^{\alpha_1},\ldots,
a^{\alpha_{p-1}}).$$ It follows from the latter equality that
$$b^{\mathcal{C}}=c.$$
\end{proof}

\begin{remark}
In many cases, it is possible to prove that these subgroups are
{\em not} conjugate in $\Gamma_{\bar{\alpha}}$.
\end{remark}

Let us establish a few other properties of $EGS$-groups. First of
all, by Theorem 2 of \cite{Ros86}, all $EGS$-groups are {\em
periodic}. Also, we have

\begin{lemma}\label{commutant_for_egs}
Let $\Gamma_{\bar{\alpha}}$ be an $EGS$-group. Then
$$\Gamma_{\bar{\alpha}}/[\Gamma_{\bar{\alpha}},\Gamma_{\bar{\alpha}}]
\cong\langle a\rangle \times\langle b\rangle \times\langle
c\rangle \cong \mathbb{Z}_p\times \mathbb{Z}_p\times
\mathbb{Z}_p.$$
\end{lemma}

The proof of this lemma can be found in \cite{congruence}, as well
as the proof of the following proposition.

\begin{proposition}\label{cost_egs}
Let $\Gamma_{\bar{\alpha}}$ be an $EGS$-group with nonsymmetric
accompanying vector. Then for any vertex $u$
$\rist_{\Gamma_{\bar{\alpha}}}(u)_u$ contains
$[\Gamma_{\bar{\alpha}},\Gamma_{\bar{\alpha}}]$.
\end{proposition}

We again stress that the proof comes from the inclusion,
\begin{equation}\label{rist_ind2}
\Phi_1([\Gamma_{\bar{\alpha}},\Gamma_{\bar{\alpha}}])\geq
\underbrace{[\Gamma_{\bar{\alpha}},\Gamma_{\bar{\alpha}}]\times
\ldots\times[\Gamma_{\bar{\alpha}},\Gamma_{\bar{\alpha}}]}_{p}.
\end{equation}

Relying on this fact, we also obtain

\begin{proposition}\label{ji_egs}
Let $\Gamma_{\bar{\alpha}}$ be an $EGS$-group with nonsymmetric
accompanying vector. Then $\Gamma_{\bar{\alpha}}$ is just
infinite.
\end{proposition}

\begin{proof}
Since $\Gamma_{\bar{\alpha}}$ is finitely generated periodic, then
so is its derived subgroup. Hence the second derived subgroup has
finite index. Then by Theorem~4 of \cite{Gr00}
$\Gamma_{\bar{\alpha}}$ is just infinite.
\end{proof}

By the same reasoning as with $G_{\bar{\alpha}}$, all
$\Gamma_{\bar{\alpha}}$ are branch.

\section{Profinite completions of $GGS$-groups}\label{ggs_completions}

The main goal of this section is to prove that the profinite
completions of $GGS$-groups coincide with their closures in $\Aut
T$ (the latter group is considered with the profinite topology).
We do that by establishing that all periodic $GGS$-groups possess
the congruence property. This is done by proving first that rigid
stabilizers are congruence subgroups (a necessary condition
always) and then showing that in fact every normal subgroup of
groups in question has a ``sufficiently large'' intersection with
some rigid stabilizer.

\subsection{Normal subgroups }

As we mentioned in Section~\ref{defin_ggs}, all $GGS$-groups are
recurrent, i.e. for any vertex $u$ $\st_{G_{\bar{\alpha}}}(u)_u$
coincides with $G_{\bar{\alpha}}$. Here we prove essentially that,
with the exception of a finite number of vertices, a similar
statement holds for any normal subgroup of
$G_{\bar{\alpha}}$.

\begin{lemma}\label{path_in}
Let $G_{\bar{\alpha}}$ be a periodic $GGS$-group. Then for any
$x\in G_{\bar{\alpha}}$ there is a vertex $u$ such that $\st_X(u)_u\ni b$,
where $X$ is the normal closure of $x$ in $G_{\bar{\alpha}}$.
\end{lemma}

\begin{proof}
Let us prove the statement by induction on length
$\partial_{G_{\bar{\alpha}}}(x)$. The base of induction is
$\partial_{G_{\bar{\alpha}}}(x)\leq 1$.

The elements of length $0$ are $a^i$, $i\in\mbox{GF}(p)$. We have
$$\Phi_1([a^i,b])=(a^{-\alpha_{p-(i-1)}},a^{-\alpha_{p-(i-2)}},\ldots,
a^{-\alpha_{p-1}},b^{-1},a^{-\alpha_1},\ldots,a^{-\alpha_{p-i}})
(a^{\alpha_1},a^{\alpha_2},\ldots,a^{\alpha_{p-1}},b)=$$
$$(a^{\alpha_1-\alpha_{p-(i-1)}},a^{\alpha_2-\alpha_{p-(i-2)}},\ldots,
b^{-1}a^{\alpha_i},a^{\alpha_{i+1}-\alpha_1},\ldots,a^{-\alpha_{p-i}}b).$$
If $\alpha_i=0$ or $\alpha_{p-i}=0$, then for either $u=i-1$ or
$u=p-i-1$, respectively, we have $\varphi_u(X)\ni b$.
If $\alpha_i\alpha_{p-i}\neq 0$, we proceed to case
$\partial_{G_{\bar{\alpha}}}(x)=1$.

If $\partial_{G_{\bar{\alpha}}}(x)=1$, then $x=a^ib^{\varepsilon}a^j$.
Conjugating by $a^i$ if needed, we can assume that $x=b^{\varepsilon}a^j$.

If $j=0$, we have the desired fact right away. Suppose $j\neq 0$. But in this
case $\varphi_{p-1}(x^p)=b^{\varepsilon}$. Since $p$ is prime, this completes
the base of induction.

Suppose now that $x$ is an element of length $m>1$, and for all
elements of smaller length the statement is already proven. If $x$
is in $\st_{G_{\bar{\alpha}}}(1)$, for any vertex $u$ of length
$1$ $\partial_{G_{\bar{\alpha}}}(\varphi_u(x))\leq
\frac{\partial_{G_{\bar{\alpha}}}(x)+1}{2}<\partial_{G_{\bar{\alpha}}}(x)$.
Since $G_{\bar{\alpha}}$ is recurrent, we can apply the inductive
assumption to any $\varphi_u(x)$.

Suppose that $x$ is not in $\st_{G_{\bar{\alpha}}}(1)$, i.e.
$x=ya^i$ for some $y\in\st_{G_{\bar{\alpha}}}(1)$ and some
$i\in\mathbb{GF}(p)$. Denote $\Phi_1(y)=(y_0,y_1,\ldots,y_{p-1})$.
Then we have
$$\Phi_1(x^p)=(y_0y_{\pi^i(0)}\ldots y_{\pi^{(p-1)i}(0)},
y_1y_{\pi^i(1)}\ldots y_{\pi^{(p-1)i}(1)},\ldots,
y_{p-1}y_{\pi^i(p-1)}\ldots y_{\pi^{(p-1)i}(p-1)}),$$
where $\pi$ denotes the (nontrivial) accompanying permutation
of $a$. As for any element of $\st_{G_{\bar{\alpha}}}(1)$,
there is equality
$$\partial_a(y_0)+\ldots+\partial_a(y_{p-1})=(\alpha_1+\ldots+\alpha_{p-1})
(\partial_b(y_0)+\ldots+\partial_b(y_{p-1}))=0.$$ Therefore, since
for any $i\in\mathbb{GF}(p)$ $a^i$ is transitive on vertices of
length $1$, $x^p$ is actually in $\st_{G_{\bar{\alpha}}}(2)$. On
the other hand, for any vertex $u$ of length $1$
$\partial_{G_{\bar{\alpha}}}(\varphi_u(x^p))\leq\partial_{G_{\bar{\alpha}}}(x)$.
Thus, for any vertex $v$ of length $2$
$\partial_{G_{\bar{\alpha}}}(\varphi_v(x^p))\leq
\frac{\partial_{G_{\bar{\alpha}}}(x)+1}{2}<\partial_{G_{\bar{\alpha}}}(x)$.
Again, since $G_{\bar{\alpha}}$ is recurrent, we can apply the
inductive assumption to any $\varphi_v(x^p)$.
\end{proof}

\begin{proposition}\label{norm_srezka}
Let $G_{\bar{\alpha}}$ be a periodic $GGS$-group. Then for any
$x\in G_{\bar{\alpha}}$ there is a vertex $u$ such that
$\st_X(u)_u=G_{\bar{\alpha}}$.
\end{proposition}

\begin{proof}
By Lemma~\ref{path_in}, there is a vertex $u$ such that
$\st_X(u)_u\ni b$. Since $G_{\bar{\alpha}}$ is recurrent and
$\st_X(u)$ is normal in $G_{\bar{\alpha}}$, we have $\st_X(u)_u\ni
b^{a^i}$ for any $i\in\mbox{GF}(p)$ as well. Thus, if $v=u(p-1)$
is a vertex adjacent to $u$, then $\st_X(v)_v\ni b$,
$a^{\alpha_1}$, $a^{\alpha_2}$, $\ldots$, $a^{\alpha_{p-1}}$.
Since among the numbers $\alpha_1$, $\alpha_2$, $\ldots$,
$\alpha_{p-1}$ there is a nontrivial one,
$\st_X(v)_v=G_{\bar{\alpha}}$.
\end{proof}

\subsection{Congruence subgroups }\label{congruence_ggs}

Consider an arbitrary accompanying vector
$\bar{\alpha}=(\alpha_1,\alpha_2,\ldots,\alpha_{p-1})$ of a
periodic $GGS$-group, i.e. one with
$\alpha_1+\alpha_2+\ldots+\alpha_{p-1}=0$. Denote by
$A_{\bar{\alpha}}$ the following matrix
\[
A_{\bar{\alpha}}=\left(
\begin{array}{cccc}
0 & \alpha_{p-1} & \ldots & \alpha_1 \\
\alpha_1 & 0 & \ldots & \alpha_2 \\
\alpha_2 & \alpha_1 & \ldots & \alpha_3 \\
\ldots & \ldots & \ldots & \ldots \\
\alpha_{p-2} & \alpha_{p-3} & \ldots & \alpha_{p-1} \\
\alpha_{p-1} & \alpha_{p-2} & \ldots & 0
\end{array}\right).
\]
This matrix can be viewed as a linear transformation of the vector
space
$\underbrace{\mathbb{GF}(p)\oplus\mathbb{GF}(p)\oplus\ldots\oplus\mathbb{GF}(p)}_{p}$.

\begin{lemma}\label{ker_of_matrix}
If a vector $(\beta_0,\ldots,\beta_{p-1})$ is in the kernel of
$A_{\bar{\alpha}}$ then the sum of its coordinates is zero,
$\sum_{i=0}^{p-1}\beta_i=0$.
\end{lemma}

\begin{proof}
Notice that the rank of $A_{\bar{\alpha}}$ is at least $1$, i.e.
the dimension of the kernel is no greater than $p-1$. Consider a
vector $(\beta_0,\ldots,\beta_{p-1})$ in the kernel. Since each
row of $A$ is obtained by a cyclic shift from any other row, the
cyclic space generated by $(\beta_0,\ldots,\beta_{p-1})$ is in the
kernel as well. On the other hand, the dimension of this cyclic
space is equal to $p-\deg\mbox{GCD}(x^p-1,f(x))$, where
$f(x)=\beta_0+\beta_1x+\ldots+ \beta_{p-1}x^{p-1}$.

Since $p$ is the characteristic of $\mathbb{GF}(p)$, there is
equality $x^p-1=(x-1)^p$. Suppose that
$\sum_{i=0}^{p-1}\beta_i\neq 0$. Then $1$ is not a root of $f(x)$,
hence $\mbox{GCD}(x^p-1,f(x))=1$. Thus, the dimension of the
cyclic space generated by $(\beta_0,\ldots,\beta_{p-1})$ is $p$.
But the dimension of the kernel, which contains this cyclic space,
does not exceed $p-1$. This contradiction proves that for any
vector from the kernel the sum of its coordinates is equal to
zero.
\end{proof}

We now easily establish the following lemma.

\begin{lemma}\label{cong_comm_ggs}
The derived subgroup of any periodic $GGS$-group $G_{\bar{\alpha}}$
contains principal congruence subgroup $\st_{G_{\bar{\alpha}}}(2)$.
\end{lemma}

\begin{proof}
Consider an element $x\in\st_{G_{\bar{\alpha}}}(1)$; denote
$\Phi_1(x)=(x_0,x_1,\ldots,x_{p-1})$. Clearly,
$\partial_b(x)=\sum_{i=0}^{p-1}\partial_b(x_i)$, and vectors
$(\partial_b(x_0),\ldots,\partial_b(x_{p-1}))$ and
$(\partial_a(x_0),\ldots,\partial_a(x_{p-1}))$ are related by
the following expression,
\[
\left(
\begin{array}{cccc}
0 & \alpha_{p-1} & \ldots & \alpha_1 \\
\alpha_1 & 0 & \ldots & \alpha_2 \\
\alpha_2 & \alpha_1 & \ldots & \alpha_3 \\
\ldots & \ldots & \ldots & \ldots \\
\alpha_{p-2} & \alpha_{p-3} & \ldots & \alpha_{p-1} \\
\alpha_{p-1} & \alpha_{p-2} & \ldots & 0
\end{array}\right)
\left(
\begin{array}{c}
\partial_b(x_0)\\
\partial_b(x_1)\\
\partial_b(x_2)\\
\ldots\\
\partial_b(x_{p-2})\\
\partial_b(x_{p-1})
\end{array}\right)=
\left(
\begin{array}{c}
\partial_a(x_0)\\
\partial_a(x_1)\\
\partial_a(x_2)\\
\ldots\\
\partial_a(x_{p-2})\\
\partial_a(x_{p-1})
\end{array}\right).\]
For $x$ to be in $\st_{G_{\bar{\alpha}}}(2)$, it is necessary and
sufficient that $\partial_a(x_i)=0$ for all $i=0,1,\ldots,p-1$.
Thus, $x\in\st_{G_{\bar{\alpha}}}(2)$ if and only if vector
$(\partial_b(x_0),\ldots,\partial_b(x_{p-1}))$ over $\mbox{GF}(p)$
is in the kernel of $A_{\bar{\alpha}}$. Then by
Lemma~\ref{ker_of_matrix} the sum of coordinates of
$(\partial_b(x_0),\ldots,\partial_b(x_{p-1}))$ is equal to zero.
Hence, $\partial_b(x)=0$, which, given that $x$ is in
$\st_{G_{\bar{\alpha}}}(1)$, implies
$x\in[G_{\bar{\alpha}},G_{\bar{\alpha}}]$.
\end{proof}

\begin{theorem}\label{cong_nonsymggs}
Any periodic $GGS$-group with a nonsymmetric accompanying vector
has congruence property.
\end{theorem}

\begin{proof}
We need to prove that every finite index subgroup of a given
$G_{\bar{\alpha}}$ is a congruence subgroup. This is equivalent to
proving that every normal finite index subgroup is a congruence
subgroup, and since all $G_{\bar{\alpha}}$ are just infinite, we
can simply prove that the normal closure $X$ of any element $x\in
G_{\bar{\alpha}}$ is a congruence subgroup.

By Proposition~\ref{norm_srezka}, there is a vertex $u$ such that
$\st_X(u)_u=G_{\bar{\alpha}}$. Notice that
$[\st_X(u),\rist_{G_{\bar{\alpha}}}(u)]\leq\rist_X(u)$.
By Corollary~\ref{cost_ggs}, $\rist_{G_{\bar{\alpha}}}(u)_u\geq
[G_{\bar{\alpha}},G_{\bar{\alpha}}]$. Hence,
$\rist_X(u)_u\geq\varphi_u([\st_X(u),\rist_{G_{\bar{\alpha}}}(u)])
\geq\gamma_3(G_{\bar{\alpha}})$. By Proposition~\ref{costt_ggs},
$\gamma_3(G_{\bar{\alpha}})\geq
\Phi_1^{-1}(\underbrace{[G_{\bar{\alpha}},G_{\bar{\alpha}}]\times\ldots
\times[G_{\bar{\alpha}},G_{\bar{\alpha}}]}_p)$. Hence, there is inclusion
$$X\geq\Phi_{|u|+1}^{-1}(\underbrace{[G_{\bar{\alpha}},G_{\bar{\alpha}}]
\times\ldots\times[G_{\bar{\alpha}},G_{\bar{\alpha}}]}_{p^{|u|+1}}).$$
However, by Lemma~\ref{cong_comm_ggs} $[G_{\bar{\alpha}},G_{\bar{\alpha}}]$
contains $\st_{G_{\bar{\alpha}}}(2)$. Therefore, $X$ contains
$\Phi_{|u|+1}^{-1}(\underbrace{\st_{G_{\bar{\alpha}}}(2)
\times\ldots\times\st_{G_{\bar{\alpha}}}(2)}_{p^{|u|+1}})=
\st_{G_{\bar{\alpha}}}(|u|+3)$. Theorem is proven.
\end{proof}

It follows from the above theorem and Proposition~2 of
\cite{GrigProf} that the closure $\bar{G}_{\bar{\alpha}}$ of
$G_{\bar{\alpha}}$ in $\Aut T$ coincides with its profinite
completion $\hat{G}_{\bar{\alpha}}$. Thus, we have the following
description of profinite completions of periodic $GGS$-groups with
nonsymmetric accompanying vectors.

\begin{theorem}\label{GGS_profinite}
The profinite completion $\hat{G}_{\bar{\alpha}}$ of a periodic
$GGS$-group with nonsymmetric accompanying vector is isomorphic to
the projective limit of the following inverse systems of finite
$p$-groups,
\begin{equation*}\begin{CD}
1 @<\pi_1<< G_{\bar{\alpha}}/\st_{G_{\bar{\alpha}}}(1) @<\pi_2<<
G_{\bar{\alpha}}/\st_{G_{\bar{\alpha}}}(2) @<\pi_3<< \ldots
@<\pi_n<<
G_{\bar{\alpha}}/\st_{G_{\bar{\alpha}}}(n) @<\pi_{n+1}<< \ldots,\\
\end{CD}\end{equation*}
where $\pi_i$'s are the natural projections.
\end{theorem}

\begin{remark}
A similar theorem can be proven for periodic $GGS$-groups with
symmetric accompanying vectors, see \cite{congruence}.
\end{remark}

\section{Completions and closures of $EGS$-groups}\label{egs_completions}

In this section we again treat $EGS$-groups with nonsymmetric
accompanying vectors only, often not specifying that.

\subsection{Absence of congruence property}

Here we establish that $EGS$-groups do {\em not} have the
congruence property.

\begin{lemma}\label{coset}
Let $\Gamma_{\bar{\alpha}}$ be an $EGS$-group with nonsymmetric
accompanying vector. Then for any natural number $n$ there is an
element $t_n$ of coset
$b[\Gamma_{\bar{\alpha}},\Gamma_{\bar{\alpha}}]$ such that
$t_n\equiv c(\mod\st_{\Gamma_{\bar{\alpha}}}(n))$.
\end{lemma}

\begin{proof}
Put $t_1=b$ and $t_2=b^a$. Suppose $n>1$ and an element $t_n$ with
the required property has already been defined. Notice that
$b^{-1}t_n$ is in $[\Gamma_{\bar{\alpha}},\Gamma_{\bar{\alpha}}]$.
Then by Proposition~\ref{cost_egs} automorphism
$x_{n+1}=\Phi_1^{-1}(b^{-1}t_n,\underbrace{1,\ldots,1}_{p-1})$ is
in $\Gamma_{\bar{\alpha}}$. Put $t_{n+1}=b^ax_{n+1}$.

We have
$$\Phi_1(t_{n+1})=(t_n,a^{\alpha_1},\ldots,a^{\alpha_{p-1}}).$$
Thus,
$$\Phi_1(c^{-1}t_{n+1})=(c^{-1}t_n,1,\ldots,1).$$
Since by assumption $c^{-1}t_n$ is in $\st_{\Gamma_{\bar{\alpha}}}(n)$,
$c^{-1}t_{n+1}$ is in $\st_{\Gamma_{\bar{\alpha}}}(n+1)$.

Finally, it follows from the proof of Proposition~\ref{cost_egs} that
$x_{n+1}$ is in $[\Gamma_{\bar{\alpha}},\Gamma_{\bar{\alpha}}]$.
Therefore $t_{n+1}$ is in $b[\Gamma_{\bar{\alpha}},\Gamma_{\bar{\alpha}}]$.
\end{proof}

\begin{corollary}\label{no_congr}
The derived subgroup $[\Gamma_{\bar{\alpha}},\Gamma_{\bar{\alpha}}]$
is not a congruence subgroup.
\end{corollary}

\begin{proof}
It follows from Lemma~\ref{conjugate} that for every $n$ $c^{-1}t_n$
is in $\st_{\Gamma_{\bar{\alpha}}}(n)$. However, by
Lemma~\ref{commutant_for_egs} it is not in
$[\Gamma_{\bar{\alpha}},\Gamma_{\bar{\alpha}}]$. Therefore,
$[\Gamma_{\bar{\alpha}},\Gamma_{\bar{\alpha}}]$ does not contain
$\st_{\Gamma_{\bar{\alpha}}}(n)$ for any $n$.
\end{proof}

Thus, we obtain the following important result.

\begin{theorem}\label{no_cong_property}
No periodic $EGS$-group with nonsymmetric accompanying vector has
congruence property.
\end{theorem}

This theorem suggests that for an $EGS$-group
$\Gamma_{\bar{\alpha}}$, its profinite completion and its closure
in $\Aut T$ are distinct. However, there is a natural homomorphism
of the former onto the latter, which is investigated in
Section~\ref{kernel}.

\subsection{Description of the profinite completions}

In this section we obtain a description of profinite completions
of $EGS$-groups as projective limits of certain linearly ordered
inverse systems of finite groups. Let us introduce a few extra
notations first. Denote by $\mathcal{H}_n$, $n\geq 1$, the
subgroup of $\Gamma_{\bar{\alpha}}$ such that
$$\Phi_n(\mathcal{H}_n)=\underbrace{[\Gamma_{\bar{\alpha}},\Gamma_{\bar{\alpha}}]
\times\ldots\times[\Gamma_{\bar{\alpha}},\Gamma_{\bar{\alpha}}]}_{p^n},$$
by $\mathcal{T}_n$ the subgroup of
$G_{\bar{\alpha}}\leq\Gamma_{\bar{\alpha}}$ given by
$$\Phi_n(\mathcal{T}_n)=\underbrace{[G_{\bar{\alpha}},G_{\bar{\alpha}}]
\times\ldots\times[G_{\bar{\alpha}},G_{\bar{\alpha}}]}_{p^n},$$
and by $\mathcal{R}_n$ the subgroup of
$F_{\bar{\alpha}}\leq\Gamma_{\bar{\alpha}}$ given by
$$\Phi_n(\mathcal{R}_n)=\underbrace{[F_{\bar{\alpha}},F_{\bar{\alpha}}]
\times\ldots\times[F_{\bar{\alpha}},F_{\bar{\alpha}}]}_{p^n}.$$ It
is also natural to denote by $\mathcal{H}_0$, $\mathcal{T}_0$, and
$\mathcal{R}_0$ subgroups
$[\Gamma_{\bar{\alpha}},\Gamma_{\bar{\alpha}}]$,
$[G_{\bar{\alpha}},G_{\bar{\alpha}}]$, and
$[F_{\bar{\alpha}},F_{\bar{\alpha}}]$, respectively.

Subgroups of $\Aut T$ $\mathcal{H}_n$, $\mathcal{T}_n$,
$\mathcal{R}_n$ are contained in groups $\Gamma_{\bar{\alpha}}$,
$G_{\bar{\alpha}}$, and $F_{\bar{\alpha}}$, respectively, by
Propositions~\ref{cost_egs}, \ref{cost_ggs} and
Lemma~\ref{conjugate}. Notice that by formulas (\ref{rist_ind2})
and (\ref{rist_ggs}) and Lemma~\ref{conjugate} these subgroups
form descending chains,
$$\mathcal{H}_0\geq \mathcal{H}_1\geq \mathcal{H}_2\geq \ldots, $$
$$\mathcal{T}_0\geq \mathcal{T}_1\geq \mathcal{T}_2\geq \ldots, $$
$$\mathcal{R}_0\geq \mathcal{R}_1\geq \mathcal{R}_2\geq \ldots. $$

The main result of this subsection is based on the following fact.

\begin{theorem}\label{small_cong_property}
Let $\Gamma_{\bar{\alpha}}$ be an $EGS$-group with nonsymmetric
accompanying vector. Then for every normal subgroup $N$ in
$\Gamma_{\bar{\alpha}}$ there is $n$ such that $N$ contains
$\mathcal{H}_n$.
\end{theorem}

The proof of this theorem is decomposed into several easy steps.
It generally follows the same scheme of reasoning as the proof of
Theorem~\ref{cong_nonsymggs}.

\begin{lemma}\label{sm_cong_g}
The normal closure in $\Gamma_{\bar{\alpha}}$ of the subgroup
$[G_{\bar{\alpha}},G_{\bar{\alpha}}]$ contains $\mathcal{H}_1$.
\end{lemma}

\begin{proof}
Notice that $[G_{\bar{\alpha}},G_{\bar{\alpha}}]$ contains all
subgroups $\mathcal{T}_n$. Choose an $i$ such that
$\alpha_1\neq\alpha_{p-i}$, which is possible to do because the
accompanying vector is not symmetric, hence not constant. Now by
direct calculation
$$[[a^{i+1},b],c]=\Phi_1^{-1}([a^{\alpha_1-\alpha_{p-i}},c],1,\ldots,1,
[b^{-1}a^{\alpha_{i+1}},a^{\alpha_{i}}],1,\ldots,1,[a^{-\alpha_{p-i-1}}b,a^{\alpha_{p-1}}])\equiv$$
$$\equiv\Phi_1^{-1}([a^{\alpha_1-\alpha_{p-i}},c],1,\ldots,1)(\mod
\mathcal{T}_1).
$$
The latter equality means that the normal closure of
$[G_{\bar{\alpha}},G_{\bar{\alpha}}]$ in $\Gamma_{\bar{\alpha}}$
contains $\mathcal{R}_1$.

Now we have
$$[[a,b],c]=\Phi_1^{-1}([b^{-1}a^{\alpha_1},c],1,\ldots,1,[a^{-\alpha_{p-1}}b,a^{\alpha_{p-1}}])\equiv
\Phi_1^{-1}([b^{-1}a^{\alpha_1},c],1,\ldots,1)(\mod
\mathcal{T}_1)\equiv
$$
$$\equiv\Phi_1^{-1}([b^{-1},c]^{a^{\alpha_1}},1,\ldots,1)(\mod
\mathcal{R}_1).
$$
Now, since $\Gamma_{\bar{\alpha}}$ is generated by $a,b,c$,
$[\Gamma_{\bar{\alpha}},\Gamma_{\bar{\alpha}}]$ is the normal
closure in $\Gamma_{\bar{\alpha}}$ of elements $[a,b]$, $[a,c]$,
$[b,c]$. Thus, since $\Gamma_{\bar{\alpha}}$ is recurrent, the two
equivalences above imply the statement of the lemma.
\end{proof}\vspace{5mm}

It is easy to obtain a similar statement for
$[F_{\bar{\alpha}},F_{\bar{\alpha}}]$ as well: 

\begin{lemma}\label{sm_cong_f}
The normal closure in $\Gamma_{\bar{\alpha}}$ of the subgroup
$[F_{\bar{\alpha}},F_{\bar{\alpha}}]$ contains $\mathcal{H}_1$.
\end{lemma}

\begin{proof}
Notice that $[F_{\bar{\alpha}},F_{\bar{\alpha}}]$ contains all
subgroups $\mathcal{R}_n$. Since the accompanying vector is
nonsymmetric and therefore non constant, there is $i$ such that
$\alpha_{i}\neq\alpha_{i+1}$. We have that
$$[[a,c],b^{a^{i+1}}]=\Phi_1^{-1}(1,\ldots,1,[a^{\alpha_{i}-\alpha_{i+1}},b],1,\ldots,1).$$
Hence the normal closure of $[F_{\bar{\alpha}},F_{\bar{\alpha}}]$
contains $\mathcal{T}_1$. On the other hand, we have
$$[[a,c],b^a]=\Phi_1^{-1}([a^{-\alpha_{p-1}}c,b],[c^{-1}a^{\alpha_1},a^{\alpha_1}],1,\ldots,1)\equiv $$
$$\equiv\Phi_1^{-1}([a^{-\alpha_{p-1}}c,b],1,\ldots,1)(\mod
\mathcal{R}_1)\equiv \Phi_1^{-1}([c,b],1,\ldots,1)(\mod
\mathcal{T}_1^c).$$ Now the desired conclusion follows at once.
\end{proof}

Combining Lemmas~\ref{sm_cong_g}, \ref{sm_cong_f}, and definitions
of subgroups $\mathcal{T}_n$, $\mathcal{R}_n$, we get the
following corollary.

\begin{corollary}\label{ryba_ili_myaso}
The normal closures in $\Gamma_{\bar{\alpha}}$ of subgroups
$\mathcal{T}_n$ and $\mathcal{R}_n$ contain $\mathcal{H}_{n+1}$.
\end{corollary}

Corollary~\ref{ryba_ili_myaso} will allow us to easily deal with
the normal subgroups contained in either $G_{\bar{\alpha}}$ or
$F_{\bar{\alpha}}$. For remaining ones, we need another
preparatory statement.

\begin{lemma}\label{ni_ryba_ni_myaso}
Let $X$ be the normal closure in $\Gamma_{\bar{\alpha}}$ of an
element $x\in\Gamma_{\bar{\alpha}}\setminus(G_{\bar{\alpha}}\cup
F_{\bar{\alpha}})$. Then there exists a vertex $u$ such that
$\st_X(u)_u$ contains at least one of $b$, $c$.
\end{lemma}

\begin{proof}
Let us prove the statement by induction on $\partial(x)$. Since
$x$ is not in $G_{\bar{\alpha}}$ or $F_{\bar{\alpha}}$, its length
is at least $2$. This is the base of the induction.

If $\partial(x)=2$ then, up to conjugation,
$x=a^{\alpha}b^{\beta}a^{\delta}c^{\gamma}$, where $\alpha$,
$\delta$ can be equal to $0$. If $\alpha+\delta=0$ then we have
two cases. One is $\alpha=-1$ and $\beta=-\gamma$. In this case
$\Phi_1(x)=(b^{-\gamma}c^{\gamma},1,\ldots,1)$, and we pass to the
second case for $x=b^{-\gamma}c^{\gamma}$.

The second case is $\alpha\neq -1$ or $\beta\neq-\gamma$. In this
case there is a vertex of length 1 such that $x_u$ is a nontrivial
element of either $G_{\bar{\alpha}}$ or $F_{\bar{\alpha}}$, and
the conclusion follows from Lemma~\ref{path_in}, or from it and
Lemma~\ref{conjugate}. If $\alpha+\delta\neq 0$ then for any $u$
of length 1 $(x^p)_u$ is, up to conjugation, either
$a^{-\mu}b^{\beta}a^{\mu}c^{\gamma}$ or
$a^{-\mu}c^{\gamma}a^{\mu}b^{\beta}$ ($\mu$ could be zero), and
the same reasoning applies.

Suppose now that for all elements of length $\leq n$ the statement
is proven, and consider $x$ of length $n+1$. If $x$ is in
$\st_{\Gamma_{\bar{\alpha}}}(1)$ then for any vertex $u$ of length
1 $x_u$ has length $\leq\frac{\partial(x)+1}{2}<\partial(x)$, and
for $x_u$ the statement is true either by the inductive assumption
(if it is not in $G_{\bar{\alpha}}\cup F_{\bar{\alpha}}$) or by
Lemmas~\ref{path_in} and \ref{conjugate} (if it is). If $x$ is not
in $\st_{\Gamma_{\bar{\alpha}}}(1)$ then, by periodicity, $x^p$
is, and $\partial(x^p)\leq
\partial(x)$. Thus, the inductive assumption or Lemmas~\ref{path_in},
\ref{conjugate} can be applied to $(x^p)_u$ for any $u$ of length
2.
\end{proof}\vspace{3mm}

\noindent\textbf{Proof of Theorem~\ref{small_cong_property} }
Obviously it is sufficient to prove the theorem for normal
closures of all elements $x\in \Gamma_{\bar{\alpha}}$. Several
cases are possible depending on whether $x$ is in
$G_{\bar{\alpha}}$, in $F_{\bar{\alpha}}$, or is not in either of
them.\vspace{3mm}

\noindent\textbf{Case $x\in G_{\bar{\alpha}}$. } Then by
Theorem~\ref{cong_nonsymggs} there is $n$ such that
$X\geq\mathcal{T}_n$. Now by Corollary~\ref{ryba_ili_myaso} $X$
contains $\mathcal{H}_{n+1}$.\vspace{3mm}

\noindent\textbf{Case $x\in F_{\bar{\alpha}}$. } Then by
Theorem~\ref{cong_nonsymggs} and Lemma~\ref{conjugate} there is
$n$ such that $X\geq\mathcal{R}_n$. By
Corollary~\ref{ryba_ili_myaso} $X$ contains
$\mathcal{H}_{n+1}$.\vspace{3mm}

\noindent\textbf{Case
$x\in\Gamma_{\bar{\alpha}}\setminus(G_{\bar{\alpha}}\cup
F_{\bar{\alpha}})$. } Then by Lemma~\ref{ni_ryba_ni_myaso} there
is a vertex $u$ such that $\st_X(u)_u$ contains either $b$ or $c$.
Suppose it contains $b$. Then it follows from the proof of
Proposition~\ref{norm_srezka} that there is a vertex $v$ (in fact,
it is one of the vertices adjacent to $u$) such that $\st_X(v)_v$
contains $G_{\bar{\alpha}}$. Applying the same reasoning as in the
proof of Theorem~\ref{cong_nonsymggs}, we get that $X$ contains
$\mathcal{T}_{|v|+1}$. Hence by Corollary~\ref{ryba_ili_myaso} $X$
contains $\mathcal{H}_{|v|+2}$. In case when $\st_X(u)_u$ contains
$c$ and not $b$, all similar arguments, with respective
substitutions, apply. \qed\vspace{5mm}

An immediate corollary of the previous theorem is the following
one.

\begin{theorem}\label{EGScompletion}
The profinite completion $\hat{\Gamma}_{\bar{\alpha}}$ is the
projective limit of the following inverse system of finite
$p$-groups:
\begin{equation*}\begin{CD}
1 @<\varepsilon_1<< \Gamma_{\bar{\alpha}}/\mathcal{H}_1
@<\varepsilon_2<< \Gamma_{\bar{\alpha}}/\mathcal{H}_2
@<\varepsilon_3<< \ldots @<\varepsilon_n<<
\Gamma_{\bar{\alpha}}/\mathcal{H}_n @<\varepsilon_{n+1}<< \ldots,\\
\end{CD}\end{equation*}
where $\varepsilon_i$ are the natural projections.
\end{theorem}

\subsection{The kernel of the natural homomorphism}\label{kernel}

\noindent\textbf{The statement of the problem and some notations }
By definition of the projective limit, the group
$\hat{\Gamma}_{\bar{\alpha}}$ is endowed with the canonical
homomorphisms $\theta_n:\hat{\Gamma}_{\bar{\alpha}}\rightarrow
\Gamma_{\bar{\alpha}}/\mathcal{H}_n$ such that
$\theta_{n-1}=\varepsilon_n\circ\theta_n$. The group
$\Gamma_{\bar{\alpha}}$ is also equipped with canonical
projections $\pi_{\mathcal{H}_n}$ onto its factor-group
$\Gamma_{\bar{\alpha}}/\mathcal{H}_n$, and the projections satisfy
the property of
$\pi_{\mathcal{H}_{n-1}}=\varepsilon_n\circ\pi_{\mathcal{H}_n}$.
Hence there is a uniquely defined homomorphism
$\hat{\sigma}:\Gamma_{\bar{\alpha}}\rightarrow
\hat{\Gamma}_{\bar{\alpha}}$ with the property that
$\pi_{\mathcal{H}_n}=\theta_n\circ\hat{\sigma}$. The kernel of
this homomorphism is the intersection
$\cap_{i=1}^{\infty}\mathcal{H}_i$. Since this intersection is
trivial, $\hat{\sigma}$ is an inclusion.

Since every finite index subgroup of $\Aut T$ contains $\st_{\Aut
T}(n)$ for some $n$, the closure $\bar{\Gamma}_{\bar{\alpha}}$ of
$\Gamma_{\bar{\alpha}}$ in $\Aut T$ considered with the profinite
topology, is the projective limit of the following inverse system
of finite $p$-groups:
\begin{equation*}\begin{CD}
1 @<\eta_1<< \Gamma_{\bar{\alpha}}/\st_{\Gamma_{\bar{\alpha}}}(1)
@<\eta_2<< \Gamma_{\bar{\alpha}}/\st_{\Gamma_{\bar{\alpha}}}(2)
@<\eta_3<< \ldots @<\eta_n<<
\Gamma_{\bar{\alpha}}/\st_{\Gamma_{\bar{\alpha}}}(n) @<\eta_{n+1}<< \ldots\mbox{ .}\\
\end{CD}\end{equation*}
Similarly, the group $\bar{\Gamma}_{\bar{\alpha}}$ is endowed with
the canonical homomorphisms
$\vartheta_n:\bar{\Gamma}_{\bar{\alpha}}\rightarrow
\Gamma_{\bar{\alpha}}/\st_{\Gamma_{\bar{\alpha}}}(n)$ such that
$\vartheta_{n-1}=\eta_n\circ\vartheta_n$. There is, too, a unique
homomorphism $\bar{\sigma}:\Gamma_{\bar{\alpha}}\rightarrow
\bar{\Gamma}_{\bar{\alpha}}$ with the property that
$\pi_{\st_{\Gamma_{\bar{\alpha}}}(n)}=\vartheta_n\circ
\bar{\sigma}$ and the trivial kernel
$\cap_{i=1}^{\infty}\st_{\Gamma_{\bar{\alpha}}}(i)$.

Finally, there are natural projections $\varphi_n:
\Gamma_{\bar{\alpha}}/\mathcal{H}_n\rightarrow
\Gamma_{\bar{\alpha}}/\st_{\Gamma_{\bar{\alpha}}}(n)$ such that
there is a commutative diagram
\begin{equation*}\begin{CD}
\Gamma_{\bar{\alpha}}/\mathcal{H}_{n-1} @<\varepsilon_n<<
\Gamma_{\bar{\alpha}}/\mathcal{H}_n \\
@VV\varphi_{n-1}V @VV\varphi_nV \\
\Gamma_{\bar{\alpha}}/\st_{\Gamma_{\bar{\alpha}}}(n-1) @<\eta_n<<
\Gamma_{\bar{\alpha}}/\st_{\Gamma_{\bar{\alpha}}}(n) \mbox{ .}
\end{CD}\end{equation*}
Therefore the set of covering homomorphisms
$\varphi_n\circ\theta_n$ allows to define a \textbf{unique
homomorphism
$\sigma:\hat{\Gamma}_{\bar{\alpha}}\rightarrow\bar{\Gamma}_{\bar{\alpha}}$
such that $\varphi_n\circ\theta_n=\vartheta_n\circ\sigma$ for all
$n$.} This homomorphism, or rather, its kernel, is the main object
of study in this section.

It is useful to remember that the three described homomorphisms
are related by the equality
$\sigma\circ\hat{\sigma}=\bar{\sigma}$. \vspace{3mm}

\noindent\textbf{One presentation of elements of the kernel } The
main goal of this section is to obtain a description of
$\Ker\sigma$. Since $\hat{\Gamma}_{\bar{\alpha}}$ is a completion
of $\Gamma_{\bar{\alpha}}$, for every
$x\in\hat{\Gamma}_{\bar{\alpha}}$ there exists a sequence
$\{g_n\}_{n=1}^{\infty}$, $g_n\in \Gamma_{\bar{\alpha}}$, such
that the sequence $\hat{\sigma}(g_n)$ is converging and its limit
is equal to $x$. Recall that $\hat{\sigma}(g_n)$ is converging if
and only if for every $N$ there is $n_N$ such that for all $n\geq
n_N$ $g_n\mathcal{H}_N=g_{n_N}\mathcal{H}_N$. Likewise, for every
$y\in\bar{\Gamma}_{\bar{\alpha}}$ there exists a sequence
$\{h_n\}_{n=1}^{\infty}$, $h_n\in \Gamma_{\bar{\alpha}}$, such
that the sequence $\bar{\sigma}(h_n)$ is converging and its limit
is equal to $y$. Here $\bar{\sigma}(h_n)$ is converging if and
only if for every $N$ there is $n_N$ such that for all $n\geq n_N$
$h_n\st_{\Gamma_{\bar{\alpha}}}(N)=h_{n_N}\st_{\Gamma_{\bar{\alpha}}}(N)$.

Thus, $x=\lim_{n\rightarrow\infty}\hat{\sigma}(g_n)$ is in
$\Ker\sigma$ if and only if the following two conditions hold:
\begin{enumerate}
\item for every $N$ there exists $n_N$ such that for all $n\geq
n_N$ $g_{n}\in g_{n_N}\mathcal{H}_N$; \item for every $N$ there
exists $m_N$ such that for all $m\geq m_N$ $g_{m}\in
\st_{\Gamma_{\bar{\alpha}}}(N)$.
\end{enumerate}
Replacing $n_N$ with $\max\{n_N,m_{N+1}\}$, we get the following
condition:
$$\forall N \mbox{ there is }n_N \mbox{ such that }
g_{n}\in g_{n_N}\mathcal{H}_N\cap
\st_{\Gamma_{\bar{\alpha}}}(N+1)\mbox{ for all }n\mbox{ greater
than }n_N.$$

Since we can assume that $n_{N-1}>n_N$, this condition implies
that $g_n\in(\bigcap_{i=1}^Ng_{n_i}\mathcal{H}_i)\cap
\st_{\Gamma_{\bar{\alpha}}}(N+1)$. Thus, we can extract a
subsequence $\{g_{n_i}\}_{i=1}^{\infty}$ such that for every $N$
$g_{n_{N+1}}\in(\bigcap_{i=1}^Ng_{n_i}\mathcal{H}_i)\cap
\st_{\Gamma_{\bar{\alpha}}}(N+1)$. Since a profinite group is
Hausdorff as a topological space, the limit of a converging
sequence coincides with the limit of every its subsequence.
Therefore $\Ker\sigma$ consists of limits of all converging
sequences $\{\hat{\sigma}(g_n)\}$, where sequence
$\{g_n\}_{n=1}^{\infty}$ is such that
\begin{equation}\label{main_condition}
g_n\in\bigcap_{i=1}^{n-1}g_i\mathcal{H}_i\cap\st_{\Gamma_{\bar{\alpha}}}(n),
\end{equation}
and we only need to consider sequences satisfying
(\ref{main_condition}). More precisely, we have the following
statement.

\begin{lemma}\label{first_char}
For every sequence $\{g_n\}$ such that
$g_n\in\cap_{i=1}^{n-1}g_i\mathcal{H}_i\cap\st_{\Gamma_{\bar{\alpha}}}(n)$
for all $n$, the limit
$\lim_{n\rightarrow\infty}\hat{\sigma}(g_n)$ exists and is in
$\Ker\sigma$. And vice versa, for every $x\in\Ker\sigma$ there
exists a sequence $\{g_n\}$ satisfying condition
$g_n\in\cap_{i=1}^{n-1}g_i\mathcal{H}_i\cap\st_{\Gamma_{\bar{\alpha}}}(n)$
and such that $x=\lim_{n\rightarrow\infty}\hat{\sigma}(g_n)$.
\end{lemma}

\begin{proof}
The second statement of the lemma has just been proven. Let us
prove the first one. Obviously, if a sequence
$\{\hat{\sigma}(g_n)\}$, where $g_n$ are of the type described in
Lemma, converges then its limit is in $\Ker\sigma$. Now the only
argument that we need to ensure the existence of limit
$\lim_{n\rightarrow\infty}\hat{\sigma}(g_n)$, is the following
evident implication: if $g_n\in g_{n-1}\mathcal{H}_{n-1}$ then
$g_n\mathcal{H}_n\subset g_{n-1}\mathcal{H}_{n-1}$. Indeed, then
$g_{n+1}\in g_n\mathcal{H}_n\subset g_{n-1}\mathcal{H}_{n-1}$, and
by induction we get that $g_{n+k+1}\in g_{n+k}\mathcal{H}_{n+k}
\subset g_{n-1}\mathcal{H}_{n-1}$ for all $n\geq 1$ and and for
all $k\geq 0$. The latter statement is just a reformulation of the
convergence condition in $\hat{\Gamma}_{\bar{\alpha}}$.
\end{proof}

We can now restrict our attention to considering only sequences
satisfying (\ref{main_condition}). Which sequences are those? is
the most important question now. For instance, an immediate
corollary of (\ref{main_condition}) is that for every sequence
$\{g_n\}$ satisfying it the coset $g_n\mathcal{H}_n$ has non-empty
intersection with all stabilizers
$\st_{\Gamma_{\bar{\alpha}}}(m)$, $m\geq n$.

\begin{definition}
We say that a coset $x\mathcal{H}_n$ {\em satisfies the kernel
condition} if it has non-empty intersection with all stabilizers
$\st_{\Gamma_{\bar{\alpha}}}(m)$, $m\geq n$.
\end{definition}

By Lemma~\ref{first_char} every element of $\Ker\sigma$ can be
presented as the limit of some converging sequence
$\{\hat{\sigma}(g_n)\}$, where all $g_n$ satisfy the kernel
condition. Therefore we need to investigate precisely which cosets
do satisfy it. \vspace{3mm}

\noindent\textbf{Cosets satisfying the kernel condition} We start
with the case $n=0$, from which all other cases easily follow.

\begin{lemma}
A coset $x[\Gamma_{\bar{\alpha}},\Gamma_{\bar{\alpha}}]$ satisfies
the kernel condition if and only if for some integer $i$
$x\equiv(c^{-1}b)^i\pmod
{[\Gamma_{\bar{\alpha}},\Gamma_{\bar{\alpha}}]}$.
\end{lemma}

\begin{proof}
It follows from Lemma~\ref{coset} that the coset $c^{-1}b
[\Gamma_{\bar{\alpha}},\Gamma_{\bar{\alpha}}]$ has non-empty
intersection with $\st_{\Gamma_{\bar{\alpha}}}(m)$ for all $m$.
Hence all cosets of the form $(c^{-1}b)^i
[\Gamma_{\bar{\alpha}},\Gamma_{\bar{\alpha}}]$ have non-empty
intersection with $\st_{\Gamma_{\bar{\alpha}}}(m)$.

On the other hand, cosets of the form
$a^{\alpha}b^{\beta}c^{\gamma}[\Gamma_{\bar{\alpha}},\Gamma_{\bar{\alpha}}]$
with $\alpha\neq 0$ obviously have empty intersection with any
stabilizer. Consider a coset of the form,
$b^{\beta}c^{\gamma}[\Gamma_{\bar{\alpha}},\Gamma_{\bar{\alpha}}]$.
An arbitrary element of this coset has the form,
$b^{\beta}c^{\gamma}g$ for
$g\in[\Gamma_{\bar{\alpha}},\Gamma_{\bar{\alpha}}]$. Denoting
$\Phi_1(g)=(g_0,\ldots,g_{p-1})$ and
$\Phi_1(b^{\beta}c^{\gamma}g)=(x_0,\ldots,x_{p-1})$, we have that
\begin{equation}\label{not_st(2)}
\left.
\begin{array}{ccccccc}
\left(
\begin{array}{c}
\partial_a(x_0) \\
\partial_a(x_1) \\
\vdots \\
\partial_a(x_{p-2}) \\
\partial_a(x_{p-1}) \\
\end{array}
\right) & = & \beta \left(
\begin{array}{c}
\alpha_1 \\
\alpha_2 \\
\vdots \\
\alpha_{p-1} \\
0
\end{array}
\right) & + & \gamma \left(
\begin{array}{c}
0 \\
\alpha_1 \\
\vdots \\
\alpha_{p-2} \\
\alpha_{p-1}
\end{array}
\right) & + & \left(
\begin{array}{c}
\partial_a(g_0) \\
\partial_a(g_1) \\
\vdots \\
\partial_a(g_{p-2}) \\
\partial_a(g_{p-1}) \\
\end{array}
\right).
\end{array}
\right.
\end{equation}
Now it is easy to show (see also the proof of
Lemma~\ref{cong_comm_ggs}) that
\[
\left.
\begin{array}{cccc}
\left(
\begin{array}{c}
\partial_a(g_0) \\
\partial_a(g_1) \\
\vdots \\
\partial_a(g_{p-2}) \\
\partial_a(g_{p-1}) \\
\end{array}
\right) & = & \left(
\begin{array}{cccc}
0 & \alpha_{p-1} & \ldots & \alpha_1 \\
\alpha_1 & 0 & \ldots & \alpha_2 \\
\alpha_2 & \alpha_1 & \ldots & \alpha_3 \\
\ldots & \ldots & \ldots & \ldots \\
\alpha_{p-2} & \alpha_{p-3} & \ldots & \alpha_{p-1} \\
\alpha_{p-1} & \alpha_{p-2} & \ldots & 0
\end{array}\right)
& \left(
\begin{array}{c}
\partial_b(g_0)+\partial_c(g_0) \\
\partial_b(g_1)+\partial_c(g_1) \\
\vdots \\
\partial_b(g_{p-2})+\partial_c(g_{p-2}) \\
\partial_b(g_{p-1})+\partial_c(g_{p-1}) \\
\end{array}
\right).
\end{array}
\right.
\]
The matrix participating in the above equality is precisely
$A_{\bar{\alpha}}$ from Section~\ref{congruence_ggs}. Notice that
\[
\left.
\begin{array}{ccc}
\beta \left(
\begin{array}{c}
\alpha_1 \\
\alpha_2 \\
\vdots \\
\alpha_{p-1} \\
0
\end{array}
\right) & = & A_{\bar{\alpha}} \left(
\begin{array}{c}
0 \\
0 \\
\vdots \\
0 \\
\beta
\end{array}
\right),
\end{array}
\right.\;\;\; \left.
\begin{array}{ccc}
\gamma \left(
\begin{array}{c}
0\\
\alpha_1 \\
\vdots \\
\alpha_{p-2} \\
\alpha_{p-1}
\end{array}
\right) & = & A_{\bar{\alpha}} \left(
\begin{array}{c}
\gamma \\
0 \\
\vdots \\
0 \\
0
\end{array}
\right),
\end{array}
\right.
\]
so the equality (\ref{not_st(2)}) can be rewritten as
\[
\left.
\begin{array}{ccc}
\left(
\begin{array}{c}
\partial_a(x_0) \\
\partial_a(x_1) \\
\vdots \\
\partial_a(x_{p-2}) \\
\partial_a(x_{p-1}) \\
\end{array}
\right) & = & A_{\bar{\alpha}} \left(
\begin{array}{c}
\gamma+\partial_b(g_0)+\partial_c(g_0) \\
\partial_b(g_1)+\partial_c(g_1) \\
\vdots \\
\partial_b(g_{p-2})+\partial_c(g_{p-2}) \\
\beta+\partial_b(g_{p-1})+\partial_c(g_{p-1}) \\
\end{array}
\right).
\end{array}
\right.
\]

It follows that if
$b^{\beta}c^{\gamma}g\in\st_{\Gamma_{\bar{\alpha}}}(2)$ then
vector $( \gamma+\partial_b(g_0)+\partial_c(g_0),
\partial_b(g_1)+\partial_c(g_1), \ldots,
\partial_b(g_{p-2})+\partial_c(g_{p-2}),\beta+\partial_b(g_{p-1})+\partial_c(g_{p-1}))$
is in the kernel of the linear transformation $A_{\bar{\alpha}}$.
By Lemma~\ref{ker_of_matrix}, for every vector from that kernel,
the sum of its coordinates should be equal to zero. On the other
hand, the sum of coordinates of our vector is
$$\beta+\gamma+(\partial_b(g_0)+\ldots+\partial_b(g_{p-1}))+(\partial_c(g_0)+\ldots+\partial_c(g_{p-1}))=
\beta+\gamma+\partial_b(g)+\partial_c(g)=\beta+\gamma,$$ because
$g$ is in $[\Gamma_{\bar{\alpha}},\Gamma_{\bar{\alpha}}]$. Thus,
an element $b^{\beta}c^{\gamma}g$ is in
$\st_{\Gamma_{\bar{\alpha}}}(2)$ if and only if $\beta+\gamma=0$.
In particular, if $\beta+\gamma\neq 0$ then the intersection of
coset
$b^{\beta}c^{\gamma}[\Gamma_{\bar{\alpha}},\Gamma_{\bar{\alpha}}]$
with $\st_{\Gamma_{\bar{\alpha}}}(2)$ is empty.
\end{proof}

Notice that for every vertex $v$ there exists a unique element
$x\in\rist_{\Gamma_{\bar{\alpha}}}(v)$ such that $x_v=c^{-1}b$.
Recall that we can denote such $x$ by $(c^{-1}b)*v$. Let
$\mathcal{CH}_n$ denote the subgroup generated by all
$(c^{-1}b)*v$ for all $v$ of length $n$,
$$\mathcal{CH}_n=\langle(c^{-1}b)*v\,:\,|v|=n \rangle.$$
Notice that all $\mathcal{CH}_n$ are abelian.\pagebreak

\begin{proposition}\label{kernel_condition}
A coset $x\mathcal{H}_n$ satisfies the kernel condition if and
only if there exists an element $y\in\mathcal{CH}_n$ such that
$x\equiv y\pmod{\mathcal{H}_n}$.
\end{proposition}

\begin{proof}
If $x\mathcal{H}_n$ satisfies the kernel condition then
$x\in\st_{\Gamma_{\bar{\alpha}}}(n)$. Consider a vertex $v$ of
length $n$. Obviously, the coset
$x_v[\Gamma_{\bar{\alpha}},\Gamma_{\bar{\alpha}}]$ must satisfy
the kernel condition. Hence for each $x_v$ there is
$y_v\in\mathcal{CH}_0=\langle c^{-1}b\rangle$ such that $x_v\equiv
y_v\pmod{[\Gamma_{\bar{\alpha}},\Gamma_{\bar{\alpha}}]}$. Thus,
$$x\equiv\prod_{|v|=n}(y_v*v)\in\mathcal{CH}_n\pmod{\mathcal{H}_n}.$$
Denoting the element on the righthand side by $y$, we get the
desired statement.
\end{proof}
\vspace{3mm}

\noindent\textbf{``Canonical'' sequences } It follows from
Proposition~\ref{kernel_condition} that we only have to consider
sequences $\{g_n\}$ such that for all $n$ and for some
$y_n\in\mathcal{CH}_n$ $g_n\equiv y_n\pmod{\mathcal{H}_n}$. More
precisely, we have the following easy fact.

\begin{lemma}\label{replacement_of_sequences}
Suppose that $\{g_n\}$, $\{y_n\}$ are two sequences such that
sequences $\{\hat{\sigma}(g_n)\}$, $\{\hat{\sigma}(y_n)\}$
converge, and for all $n$ $g_n\equiv y_n\pmod{\mathcal{H}_n}$.
Then
$$\lim_{n\rightarrow\infty}\hat{\sigma}(g_n)=\lim_{n\rightarrow\infty}\hat{\sigma}(y_n).$$
\end{lemma}

\begin{proof}
The condition $g_n\equiv y_n\pmod{\mathcal{H}_n}$ for all $n$
means that, for every $N$ and for all $n\geq N$
$g_n^{-1}y_n\in\mathcal{H}_n\subset\mathcal{H}_N$. This means that
$\lim_{n\rightarrow\infty}\hat{\sigma}(g_n^{-1}y_n)$ is the
trivial element of $\hat{\Gamma}_{\bar{\alpha}}$. Since sequences
$\{\hat{\sigma}(g_n)\}$, $\{\hat{\sigma}(y_n)\}$ converge, this is
equivalent to the statement of the lemma.
\end{proof}

Hence we are only interested in sequences $\{g_n\}$ such that for
all $n$ $g_n\in\mathcal{CH}_n$. More precisely,

\begin{corollary}\label{form_of_sequences}
For every element $x$ of $\Ker\sigma$ there exists a sequence
$\{g_n\}$ such that $g_n\in\mathcal{CH}_n$ for all $n$, sequence
$\{\hat{\sigma}(g_n)\}$ converges, and
$x=\lim_{n\rightarrow\infty}\hat{\sigma}(g_n)$. Vice versa, if a
sequence $\{\hat{\sigma}(g_n)\}$ such that $g_n\in\mathcal{CH}_n$
for all $n$, converges, then its limit is in $\Ker\sigma$.
\end{corollary}

The proof of the corollary follows immediately from
Lemma~\ref{first_char}, Proposition~\ref{kernel_condition}, and
the definition of the kernel condition. Due to the corollary, we
now only need to establish when a sequence of the type described
in Corollary~\ref{form_of_sequences} converges, and when two such
sequences have the same limit. \vspace{3mm}

\noindent\textbf{Convergence of canonical sequences } Although
that is not a necessary condition of convergence, due to
Lemma~\ref{first_char} we can confine ourselves to considering
only sequences satisfying (\ref{main_condition}). Since all
$\mathcal{CH}_n$ are in $\st_{\Gamma_{\bar{\alpha}}}(n)$ already,
we just need to ensure that for all $n$ $g_n$ be in
$\cap_{i=1}^{n-1}g_i\mathcal{H}_i$. In particular, we want this
intersection to be non-empty. Recall also the previously mentioned
(Lemma~\ref{first_char}) fact that
\begin{equation}\label{cosets_and_elements}
g_n\in\bigcap_{i=1}^{n-1}g_i\mathcal{H}_i\Longleftrightarrow
g_n\mathcal{H}_n\subset \bigcap_{i=1}^{n-1}g_i\mathcal{H}_i.
\end{equation}

\begin{lemma}\label{base_of_convergence}
Let $u$, $v$ be vertices of $T$ such that $u\leq v$ and
$n=|u|=|v|-1$. Then $((c^{-1}b)*u)\mathcal{H}_n\supset
((c^{-1}b)*v)\mathcal{H}_{n+1}$.
\end{lemma}

\begin{proof}
We have, for a suitably chosen $i$,
$$((c^{-1}b)*v)_u=(c^{-1}b[b,a])^{a^i}\equiv c^{-1}b\pmod{[\Gamma_{\bar{\alpha}},\Gamma_{\bar{\alpha}}]}.$$
This means that there is
$x\in[\Gamma_{\bar{\alpha}},\Gamma_{\bar{\alpha}}]$ such that
$((c^{-1}b)*v)_u=c^{-1}bx$. Hence
$$(c^{-1}b)*v=((c^{-1}b)*u)(x*u),$$ and
since $x*u\in\mathcal{H}_n$, we have that
$$((c^{-1}b)*u)\mathcal{H}_n=((c^{-1}b)*v)\mathcal{H}_n\supset
((c^{-1}b)*v)\mathcal{H}_{n+1}.$$ This proves the lemma.
\end{proof}

\begin{corollary}\label{description}
Let $\{g_n\}$ be a sequence such that for all $n$
$g_n\in\mathcal{CH}_n$, $g_n=\prod_{|v|=n}(c^{-1}b)^{i_v}*v$. Then
the sequence $\{\hat{\sigma}(g_n)\}$ satisfies
(\ref{main_condition}) if and only if the following condition
holds: for every vertex $v$ $i_v\equiv i_{v_1}+\ldots+i_{v_p}\pmod
p$, where $v_1$, $\ldots$, $v_p$ are all vertices of length
$|v|+1$ joined to $v$ by an edge.
\end{corollary}

\begin{proof}
Consider an arbitrary $m\geq 1$, and for each vertex $v$ of length
$m-1$ denote by $v_1$, $\ldots$, $v_p$ the $p$ vertices of length
$m$ adjacent to it. In particular, each vertex of length $m$ is
written as $v_i$ for some $v$, $i$. Then by definition of numbers
$i_v$ and Lemma~\ref{base_of_convergence} we have:
$$g_m\mathcal{H}_m=\prod_{|v|=m-1}((c^{-1}b)^{i_{v_1}}*v_1\cdot\ldots(c^{-1}b)^{i_{v_p}}*v_p)\mathcal{H}_m=$$
$$=\prod_{|v|=m-1}((c^{-1}b)^{i_{v_1}}*v_1)\mathcal{H}_m\cdot\ldots((c^{-1}b)^{i_{v_p}}*v_p)\mathcal{H}_m
\subset$$
$$\subset\prod_{|v|=m-1}((c^{-1}b)^{i_{v_1}}*v)\mathcal{H}_{m-1}\cdot\ldots((c^{-1}b)^{i_{v_p}}*v)\mathcal{H}_{m-1}=$$
$$=\prod_{|v|=m-1}((c^{-1}b)^{i_{v_1}}*v\cdot\ldots((c^{-1}b)^{i_{v_p}}*v)\mathcal{H}_{m-1}=
\prod_{|v|=m-1}((c^{-1}b)^{i_{v_1}+\ldots+i_{v_p}}*v)\mathcal{H}_{m-1}.
$$

By Lemma~\ref{commutant_for_egs} cosets
$(c^{-1}b)^i[\Gamma_{\bar{\alpha}},\Gamma_{\bar{\alpha}}]$ and
$(c^{-1}b)^j[\Gamma_{\bar{\alpha}},\Gamma_{\bar{\alpha}}]$ have
non-empty intersection if and only if $i\equiv j\mod p$. Therefore
the coset
$\prod_{|v|=m-1}((c^{-1}b)^{i_{v_1}+\ldots+i_{v_p}}*v)\mathcal{H}_{m-1}$
has non-empty intersection with the coset
$g_{m-1}\mathcal{H}_{m-1}$ if and only if $v_1+\ldots+v_p\equiv
i_v\pmod p$ for all $v$. Since two cosets of the same subgroup
have non-empty intersection if and only if they coincide, we have
inclusion $g_m\mathcal{H}_m\subset g_{m-1}\mathcal{H}_{m-1}$. This
guarantees us that, for all $n$, $g_n\mathcal{H}_n$ is in
$\bigcap_{i=1}^{n-1}g_i\mathcal{H}_i$. By
(\ref{cosets_and_elements}), this is equivalent to $g_n$ being in
$\bigcap_{i=1}^{n-1}g_i\mathcal{H}_i$. Finally, since
$g_n\in\st_{\Gamma_{\bar{\alpha}}}(n)$ by definition, $g_n$
satisfies condition (\ref{main_condition}).
\end{proof}
\vspace{3mm}

\noindent\textbf{Uniqueness of canonical sequences } Any sequence
satisfying conditions of Corollary~\ref{description}, is given by
the set of indices $i_v\in\{0,1,\ldots,p-1\}$ placed at each
vertex of tree $T$ and satisfying the ``summation condition''
stated in that corollary. The next natural question to ask is, is
the collection of indices uniquely determined by the limit of the
sequence?

\begin{proposition}\label{unique_sequence}
Let $\{g_n\}$, $\{h_n\}$ be two sequences satisfying the
conditions of Corollary~\ref{description},
$g_n=\prod_{|v|=n}(c^{-1}b)^{i_v}*v$,
$h_n=\prod_{|v|=n}(c^{-1}b)^{j_v}*v$. Then
$\lim_{n\rightarrow\infty}\hat{\sigma}(g_n)=\lim_{n\rightarrow\infty}\hat{\sigma}(h_n)$
if and only if for all $v$ $i_v\equiv j_v\pmod p$.
\end{proposition}

\begin{proof}
If for all $v$ $i_v\equiv j_v\pmod p$ then
$\hat{\sigma}(g_n)=\hat{\sigma}(h_n)$, and the equality of limits
is evident. Suppose that the limits are equal. This means that
$\lim_{n\rightarrow\infty}\hat{\sigma}(g_n^{-1}h_n)$ is the
trivial element of $\hat{\Gamma}_{\bar{\alpha}}$, i.e. that for
every $N$ there is $n_N$ such that for all $n\geq n_N$
$g_n^{-1}h_n\in\mathcal{H}_N$. Consider some $n\geq\max\{n_N,N\}$. 
It was shown in the proof of Corollary~\ref{description} that
$g_n\mathcal{H}_n\subset g_N\mathcal{H}_N$,
$h_n\mathcal{H}_n\subset h_N\mathcal{H}_N$. In particular, $g_n\in
g_N\mathcal{H}_N$ and $h_n\in h_N\mathcal{H}_N$. Hence
$g_n^{-1}h_n\in (g_N^{-1}h_N)\mathcal{H}_N$, and so
$$g_n^{-1}h_n\in (g_N^{-1}h_N)\mathcal{H}_N\bigcap\mathcal{H}_N.$$
Since the latter two sets are cosets of the same subgroup, they
can only have non-empty intersection if they coincide. This means
that $g_N^{-1}h_N\in\mathcal{H}_N$, i.e. that for every vertex $u$
of length $N$ $i_u\equiv j_u\pmod p$. Since $N$ is any, this is
actually true for any vertex at all.
\end{proof}
\vspace{3mm}

\noindent\textbf{At last, description of the kernel } It follows
from Corollary~\ref{description} and
Proposition~\ref{unique_sequence} that $\Ker\sigma$ can be
described in the following way. Consider a collection of groups
$\tilde{\mathbb{Z}}_p^{(n)}=\bigoplus_{|v|=n}\mathbb{Z}_p^{(v)}$,
the $n$th of which is the direct sum of $p^n$ exemplars of the
cyclic group $\mathbb{Z}_p$. The $0$th group
$\tilde{\mathbb{Z}}_p^{(0)}$ is trivial. The exemplars are
parameterized by vertices of length $n$ in tree $T$. Let $\pi_u$
be the projection of the corresponding group
$\tilde{\mathbb{Z}}_p^{(n)}$, $n=|u|$, onto its summand
$\mathbb{Z}_p^{(u)}$. Consider the map
$\theta_n:\tilde{\mathbb{Z}}_p^{(n)}\rightarrow
\tilde{\mathbb{Z}}_p^{(n-1)}$, $n\geq 2$, defined by the rule
$$\pi_u(\theta_n(x))=\sum_{v\geq u,\,|v|=n}\pi_v(x),$$
where $u$ is any vertex of length $n-1$. The map
$\theta_1:\mathbb{Z}_p\rightarrow \{0\}$ simply sends everything
to zero.

\begin{theorem}\label{final_result}
$\Ker\sigma$ is isomorphic to the projective limit
$\tilde{\mathbb{Z}}_p^{(\infty)}$ of the following inverse system
of groups,
\begin{equation*}\begin{CD}
1 @<\theta_1<< \tilde{\mathbb{Z}}_p^{(1)} @<\theta_2<<
\tilde{\mathbb{Z}}_p^{(2)} @<\theta_3<< \ldots @<\theta_n<<
\tilde{\mathbb{Z}}_p^{(n)} @<\theta_{n+1}<< \ldots\mbox{ .}\\
\end{CD}\end{equation*}
In particular, $\Ker\sigma$ is a profinite abelian group of
(prime) exponent $p$.
\end{theorem}

\begin{proof}
We establish a map $\Lambda:\Ker\sigma\rightarrow
\tilde{\mathbb{Z}}_p^{(\infty)}$ in the following way. Let $x$ be
in $\Ker\sigma$. Then
$x=\lim_{n\rightarrow\infty}\hat{\sigma}(g_n)$ with $g_n$ as in
Corollary~\ref{description}, $g_n=\prod_{|v|=n}(c^{-1}b)^{i_v}*v$.
Consider map $\tau_n:\Ker\sigma\rightarrow
\tilde{\mathbb{Z}}_p^{(n)}$ defined by
$$\pi_v(\tau_n(x))=i_v \mbox{ for each vertex } v \mbox{ of length } n.$$
(We assume that all $i_v$ are taken from the set
$\{0,1,\ldots,p-1\}$.) Then by definition of the projective limit,
the collection of maps $\tau_n$ defines a unique map
$\Lambda:\Ker\sigma\rightarrow \tilde{\mathbb{Z}}_p^{(\infty)}$.
However, we need to make sure that all maps $\tau_n$ are
well-defined (i.e. that $i_v$ depend on $x$ only), that they are
homomorphic and that $\Lambda$ is bijective. The first two
conditions follow immediately from
Proposition~\ref{unique_sequence}. To get the bijectivity
condition, we construct an inverse map $\Omega$ by putting for
each $y\in\tilde{\mathbb{Z}}_p^{(\infty)}$
$$\Omega(y)=\lim_{n\rightarrow\infty}\hat{\sigma}(g_n),\mbox{ where
} g_n=\prod_{|v|=n}(c^{-1}b)^{\pi_v(\zeta_n(y))}*v
$$
($\zeta_n$ are the canonical maps
$\tilde{\mathbb{Z}}_p^{(\infty)}\rightarrow
\tilde{\mathbb{Z}}_p^{(n)}$). Evidently, this is a well-defined
inverse map and a homomorphism.
\end{proof}

It is interesting to note that $\Ker\sigma$ does not depend on a
particular vector $\bar{\alpha}$, only on the number $p$ (as soon
as the vector is nonsymmetric and the sum of its coordinates is
zero).

Notice also that the group $\tilde{\mathbb{Z}}_p^{(\infty)}$ is
rather large. Indeed, for every infinite (strictly descending)
path $\gamma$ in tree $T$ it contains an element $b_{\gamma}$
defined by the rule
$$\pi_v(\zeta_n(b_{\gamma}))=\left\{
\begin{array}{l}
1,\mbox{ if }v\in\gamma,\\
0,\mbox{ otherwise.}
\end{array}\right.$$
Since for every finite set of paths there is a vertex lying on
exactly one of those paths, and $p$ is prime, the subgroup
$B^{\omega}$ generated by all $b_{\gamma}$ is generated freely by
them. Hence this subgroup is isomorphic to the additive group of
the space of all finite-support functions on $\partial T$ taking
values in $\mathbb{GF}(p)$. Notice also that for any $z\in
\tilde{\mathbb{Z}}_p^{(\infty)}$
$$z\equiv\sum_{|v|=n}b_{\gamma_v}^{\pi_v(\zeta_n(z))}\pmod{\Ker\zeta_n},$$
where $\gamma_v$ is the infinite path of the form, $v00\ldots
0\ldots$. Hence $B^{\omega}$ is dense in
$\tilde{\mathbb{Z}}_p^{(\infty)}$.


\end{document}